# Optimizing Planning Service Territories by Dividing Into Compact Several Sub-areas Using Binary K-means Clustering According Vehicle Constraints

**Muhammad Wildan Abdul Hakim, Syarifah Rosita Dewi , Yurio Windiatmoko,
Umar Abdul Aziz**
hakim@warungpintar.co, syarifah@warungpintar.co, yurio@warungpintar.co, umar@warungpintar.co
Warung Pintar Data Research, Jakarta, Indonesia

**Abstract**. VRP (Vehicle Routing Problem) is an NP hard problem, and it has attracted a lot of research interest. In contexts where vehicles have limited carrying capacity, such as volume and weight but needed to deliver items at various locations. Initially before creating a route, each vehicle needs a group of delivery points that are not exceeding their maximum capacity. Drivers tend to deliver only to certain areas. Cluster-based is one of the approaches to give a basis for generating tighter routes. In this paper we propose new algorithms for producing such clusters/groups that do not exceed vehicles maximum capacity. Our basic assumptions are each vehicle originates from a depot, delivers the items to the customers and returns to the depot, also the vehicles are homogeneous. This methods are able to compact sub-areas in each cluster. Computational results demonstrate the effectiveness of our new procedures, which are able to assist users to plan service territories and vehicle routes more efficiently.

## 1. Introduction

In practical applications of routing and distribution, a service or logistics provider continually faces the challenge of providing the best service to customers. Typically this challenge is cast in the framework of determining the most effective way to deliver services from customers to other customers within a specified area, subject to limitations on resources such as service personnel and vehicle volume. Most of these algorithms employ a cluster-first and route-second strategy in order to address real application problems effectively.

Another commonly encountered challenge in solving these problems arises from the fact that, due to the limited availability of resources, a service or logistics provider is generally unable to service all customers within a service territory. The problem then arises: 'how to create these sub-areas to meet the business logic'. The service resource is limited. We are required to develop a decision-support system to help the user create sub-areas within this territory so that each sub-area can be serviced directly. Each sub-area should be as compact as possible so that the expected service cost and especially the travel time and distance to service the sub-areas will be minimal. In the following discussion we will use the term cluster to represent a sub-area.

Expressed more generally, the problem consists in determining how to group or cluster the customers, which may abstractly be treated as points within the service territory, while honoring the rules imposed by business practice such as restrictions on total service capacity (total working hours combined with vehicle capacities including weights and volumes), balanced workload (service levels or delivery quantities), non-overlapped subareas, etc.

## 2. Literature Review

In context of optimizing planning service territories into sub-areas, Dondo & Cerda (2007)[1] [1] proposed a cluster-based approach to give a basis for generating tighter routes. In this procedure, all customer locations are clustered and the clusters in turn are assigned to vehicles. The approach accounts for vehicle capacities, time windows and idle time. The goal is to build a cluster yielding a low average travel distance for each location.

Several assignment and heuristic algorithms have been developed in the past decade to solve CCP (The Capacitated Clustering Problem). Geetha, Poonthalir & Vanath (2009) [5] applies k-means algorithm to solve CCP. The Capacitated Clustering Problem (CCP) partitions a set of n items (eg. customer orders) into k disjoint clusters with known capacity. During clustering the items with shortest assigning paths from centroids are grouped together. The summation of grouped items should not exceed the capacity of the cluster. All clusters have uniform capacity. In the computational study, the benchmark data instances are solved to evaluate the tightness measure. Experimental results have shown that the priority applied to handle the capacity constraint can guide the search direction.

To improve delivery efficiency, Hosoda J. & Irohara T. (2019) [6] proposed a new location routing problem (LRP) model with outsourced delivery using clustering algorithm. The model includes not only selection of depot location and vehicle routes, but also selection of delivery mode, that is selection of non-outsourced delivery or outsourced delivery. The algorithm divides the original LRP into smaller problems to make it easier to solve, without reducing the optimization potential as much as possible. The idea of clustering is not to divide all customers into clusters, but to cluster only close customers. In the next step, non-clustered customers are assigned to vehicles. The computer experiments on small-scale problems showed that the case in which about 80% of the customers were grouped into a cluster could get better solutions.

Our methodology undertakes to provide the following contributions to handling service territory planning and design issues in the logistics and service industry. A new framework that enables users in the logistics and service industry to plan and design their service areas more efficiently because of the non computationally expensive process.

## 3. Data and Methods

Optimizing planning service territories into several sub-areas becomes our domain research. Hence, we firstly collect a simulation dataset of sales order list within the delivery points each. The attributes used in this dataset are, *origin* as the unique sale order, *vol_cbm* as the total volume of each product ordered in the sale order, *weight_ton* as the total tonnage of each product ordered on the sale order, *partner_longitude* as the longitude of delivery point, *partner_latitude* as the latitude of delivery point.

The goal of our problem, roughly stated, is to create a clustering set $\Omega$ so that each cluster C $\in \Omega$ contains approximately the same total volume cubic meters. According to a specified level of priority, create a clustering set $\Omega$ in relation to a distance measure so that each cluster $C \in \Omega$ is composed of nodes that are closer to other nodes of the same cluster than they are to nodes of other clusters (Node here's denoted as delivery point). The motivation underlying the problem objective is to produce clusters that give a foundation for creating routes that can be served by different vehicles on the same days.

---
[1] "A cluster-based optimization approach for the multi-depot ...." 1 Feb. 2007, https://www.sciencedirect.com/science/article/pii/S0377221705008672. Accessed 5 Oct. 2020.

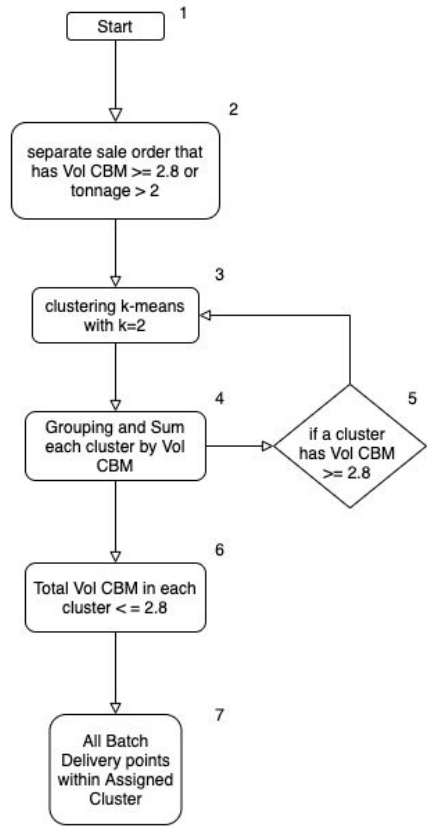

**Figure 1.** Flow Diagram Cluster by Constraint

Initially, the dataset must be separated from the single sale order that has total CBM volume >= 2.8 or with total tonnage > 2, then those sales order collection is eliminated. After that, the remaining dataset is carried out by a k-means clustering process with a value of k = 2, then the two clusters are grouped to calculate the total CBM volume, then if the volume of cbm in the cluster still exceeds or more than 2.8, the cluster will be carried out by the k-means clustering process again with a value of k = 2 and so on, until those of cluster collection has CBM volumes <= 2.8 each, this iteration process is carried out with a binary tree algorithm. The iteration with the binary tree k-means algorithm can be illustrated in the image below.

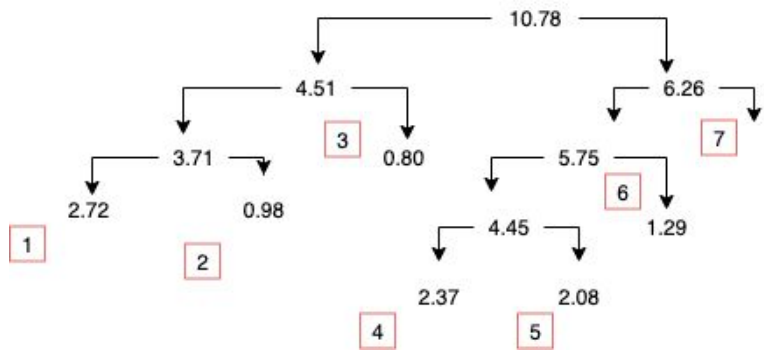

**Figure 2.** Binary Tree K-Means Illustration

Figure 2. illustrates the binary tree k-means iteration process with k = 2, which the process repeats the k-means clustering process in each cluster, where the grouping results of the volume of cbm are still more than 2.8. After that traverse all the cluster leaves (1, 2, 3, 4, 5, 6, 7 pointed aboves) to get all cluster tables, then finally clustered data already served.

### 3.1 K-Means Clustering

The K-Means algorithm clusters data by trying to separate samples in n groups of equal variance, minimizing a criterion known as the inertia or within-cluster sum-of-squares (see formulas below). This algorithm requires the number of clusters to be specified. It scales well to a large number of samples and has been used across a large range of application areas in many different fields. The k-means algorithm divides a set of samples into disjoint clusters, each described by the mean of the samples in the cluster. The means are commonly called the cluster "centroids", note that they are not in general points from, although they live in the same space. The K-means algorithm aims to choose centroids that minimise the inertia, or within-cluster sum-of-squares criterion:

$$\sum_{i=0}^{n} \min_{\mu_j \in C}(||x_i - \mu_j||^2)$$

which $(||x_i - \mu_j||^2)$ as Squared Euclidean norm of data point from *i = 0, 1, 2, . . n* and *C as a cluster*, $x_i$ *as data point x at i and* $\mu_j$ *as centroid of* $x_i$ *cluster*. We set 0 as the *random_state* of the generator to be used to initialize the centers and 2 for *n_cluster* or *k* as total cluster.[2] [3] [2] [3]

### 3.2 Custom Binary Tree K-Means Clustering

The k-means optimization problem here is to find the set of clusters which is compact and composed of points that are closer to other points of the same cluster than they are to points of other clusters also satisfying the vehicle volume constraints. We propose the use of a binary tree for k-means clustering, given in Algorithm 1 and traversing all cluster data from its leaves in Algorithm 2.[4] [4]. The motivation behind this method is to create a partition by area and satisfy the vehicle volume constraints, so that each cluster would have a tighter route. Binary tree has proven its use that does not suffer increased computational cost when datasets grow large with redundant examples. [4] [7]

---

[2] "k-means++: The Advantages of Careful Seeding - Stanford CS ...."
https://theory.stanford.edu/~sergei/papers/kMeansPP-soda.pdf. Accessed 6 Oct. 2020.
[3] "Scikit-learn: Machine Learning in Python - Journal of Machine ...."
https://jmlr.csail.mit.edu/papers/volume12/pedregosa11a/pedregosa11a.pdf. Accessed 6 Oct. 2020.
[4] "Teaching binary tree algorithms through visual ... - IEEE Xplore."
https://ieeexplore.ieee.org/document/545265. Accessed 6 Oct. 2020.

**Algorithm 1** BT-Insert

```
1: DEFINE CLASS NODE:
2:       DEFINE FUNCTION INIT DATA:
3:            SET LEFT TO NONE
4:            SET RIGHT TO NONE
5:            SET DATA TO DATA
6:            IF LENGTH(DATA) > 1:
7:                SET BRANCH, VOL_CBM_BRANCH, WEIGHT_BRANCH TO
                     CLUSTER_FUNC(DATA)
8:                INSERT_FUNC(BRANCH)
9:            ELSE:
10:               PASS
11:      DEFINE FUNCTION INSERT_FUNC(BRANCH):
12: IF MAX(VOL_CBM_BRANCH_CLUSTER_0,VOL_CBM_BRANCH_CLUSTER_1) <= 2.8:
13:          PASS
14: ELSE:
15:      IF VOL_CBM_BRANCH_CLUSTER_0 <= 2.8:
16:          PASS
17:      ELSE:
18:          SET LEFT TO NODE(BRANCH_CLUSTER_0)
19:      IF VOL_CBM_BRANCH_CLUSTER_1 <= 2.8:
20:          PASS
21:      ELSE:
22:          SET RIGHT TO NODE(BRANCH_CLUSTER_1)
```

**Algorithm 2** Traverse-Node

```
1: DEFINE FUNCTION_TRAVERSE(NODE):
2:       SET GLOBAL COUNTER
3:       SET GLOBAL LIST_CONTAINER
4:       SET GLOBAL TABLE_CONTAINER
5:       TRY:
6:            IF NODE_VOL_CBM_BRANCH_CLUSTER_0 <= 2.8:
```

```
7:                         COUNTER+=1
8:                         LIST_CONTAINER.ADD(COUNTER)
9:                         TABLE_CONTAINER.ADD(NODE_BRANCH_CLUSTER_0)
10:                 IF NODE_VOL_CBM_BRANCH_CLUSTER_1 <= 2.8:
11:                         COUNTER+=1
12:                         LIST_CONTAINER.ADD(COUNTER)
13:                         TABLE_CONTAINER.ADD(NODE_BRANCH_CLUSTER_1)
14:         EXCEPT :
15:                 PASS
16:     IF NODE IS NONE:
17:             RETURN 0
18:     IF(NODE_LEFT IS NONE AND NODE_RIGHT IS NONE):
19:             RETURN 1
20:     ELSE:
21:             RETURN FUNCTION_TRAVERSE(NODE_LEFT) +
                    FUNCTION_TRAVERSE(NODE_RIGHT)
```

## 4. Result and Discussion

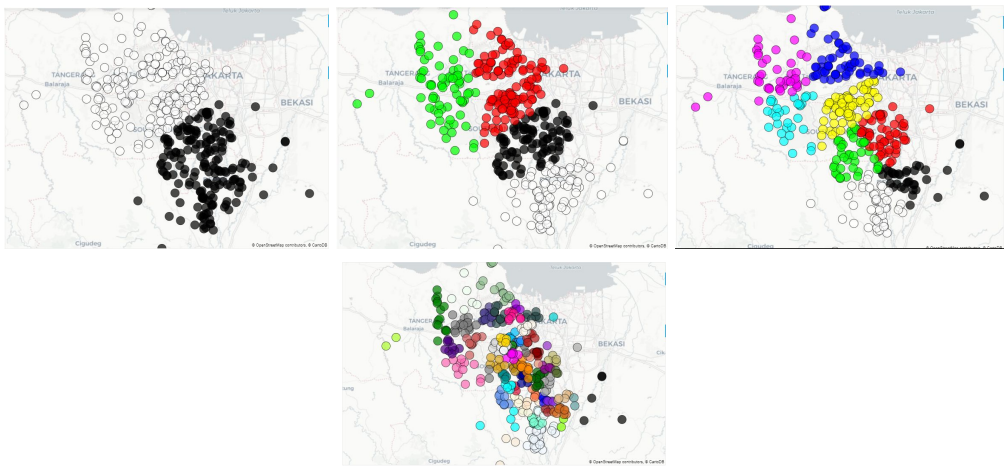

**Figure 3.** Result Plot Step by Step Cluster

Total clusters that were formed are relatively a lot and have various volumes. Indeed it has already satisfied the volume constraints but in contrast, there are a lot of clusters formed with a sum of volume very low that causes vehicles not well utilized when there is so much vehicle space left.

## 5. Conclusion and Future Work

To optimize vehicle capacity used, the future work would be to set a threshold minimum capacity assigned to vehicles after clustering. Create a binary tree cluster for another cluster below threshold. We assume it would cause a better utilization for each vehicle and also decrease mileage of the vehicles that cut delivery time.